\documentclass[11pt,reqno]{amsart} 

\usepackage{mypack}

\usepackage{latexsym}
\usepackage{a4wide}
\usepackage{amscd}
\usepackage{graphics}
\usepackage{amsmath}
\usepackage{amssymb}
\usepackage{mathrsfs}
\usepackage{accents}
\usepackage{enumerate}
\usepackage[T1]{fontenc}

\usepackage[backgroundcolor=white,linecolor=gray]{todonotes}

\DeclareMathOperator{\Th}{Th}

\newcommand{\frc}{\Omega^{\mathrm{fr}}}

\author{Thomas O. Rot}
\begin{document}


\title[Homotopy classes of proper maps]{Homotopy classes of proper maps out of vector bundles.}
\begin{abstract}
In this paper we classify the homotopy classes of proper maps $E\rightarrow \mR^k$, where $E$ is a vector bundle over a compact Hausdorff space. As a corollary we compute the homotopy classes of proper maps $\mR^n\rightarrow \mR^k$. We find a stability range of such maps. We conclude with some remarks on framed submanifolds of non-compact manifolds, the relationship with proper homotopy classes of maps and the Pontryagin-Thom construction. 
\end{abstract}
\maketitle 

\section{Introduction}
A continuous map $f:X\rightarrow Y$ is called proper if $f^{-1}(C)$ is compact for all compact subsets $C$ of $Y$. A homotopy of proper maps is a homotopy $F:[0,1]\times X\rightarrow Y$ such that $F$ is a proper map. The assumption that a homotopy is a homotopy of proper maps is stronger than the assumption that the homotopy is homotopy through proper maps, i.e. the assumption that the maps $F_t:X\rightarrow Y$ are proper for every $t\in[0,1]$. A simple example of a homotopy through proper maps that is not a homotopy of proper maps is the map $F:[0,1]\times \mR\rightarrow \mR$ defined by $F(t,x)=(1-t)x^2+x$. To see this note that the sequence $(1-\frac{1}{n},-n)$ is unbounded, but $F(1-\frac{1}{n},-n)=0$. This example is closely related to the compactness issues discussed in~\cite{Geiges:2018fo}.

We denote by $[X,Y]$ the set of (unbased) homotopy classes of maps from $X$ to $Y$ and with $[X,Y]_{\mathrm{prop}}$ the set of (unbased) homotopy classes of proper maps. For the set of homotopy classes of based maps between pointed spaces we write $\langle X,Y\rangle$. 

In~\cite{Abbondandolo:2018fb} we classified the homotopy classes of proper Fredholm maps of Hilbert manifolds into its model (real and separable) Hilbert space in terms of a suitable notion of framed cobordism. This classification uses an infinite-dimensional and proper analogue of the Pontryagin-Thom collapse map, which is due to Elworthy and Tromba~\cite{Elworthy:1970vp}, see also the paper of G\k{e}ba~\cite{Geba:1969df}. The existence of the collapse map hinges on the fact that an infinite dimensional Hilbert space is diffeomorphic to the Hilbert space minus a point. This is of course not true for a finite dimensional vector space and the construction does not work in this setting. As we will discuss in Section~\ref{framed}, even though the framed cobordism class of a regular value is an invariant of the homotopy class of a proper map in the finite dimensional setting, the framed cobordism class is not able to distinguish all proper homotopy classes of proper maps into $\mathbb{R}^k$, nor do all framed submanifolds come from proper maps. Thus there does not exist a finite dimensional proper Pontryagin-Thom construction, which is why we are not able to compute $[E,\mR^k]_{\mathrm{prop}}$ for all open finite dimensional manifolds $E$ using a Pontryagin-Thom collapse map. In this paper we are content with the classification of $[E,\mR^k]_{\mathrm{prop}}$ where $E$ is a real vector bundle over a compact Hausdorff space $M$. This classification does not use a Pontryagin-Thom collapse map.  

\begin{theorem}
\label{maintheorem}
Let $E\rightarrow M$ be a normed vector bundle over a compact Hausdorff space $M$. Denote by $S(E)\rightarrow M$ the associated sphere bundle of unit vectors. Then the set $[E,\mathbb{R}^k]_{\mathrm{prop}}$ is in bijective correspondence with $[S(E),S^{k-1}]$. 
\end{theorem}

We have the following corollary of Theorem~\ref{maintheorem} by taking $M$ to be a point and using the fact that based and unbased homotopy classes of maps from spaces to positive dimensional spheres coincide, cf.~\cite[Section 4A]{Hatcher_Topology}. 

\begin{corollary}
\label{cor:sphere}
$[\mR^n,\mR^k]_{\mathrm{prop}}$ is in bijective correspondence with $[S^{n-1},S^{k-1}]$. Thus if $n>1$ and $k>1$ we have that $[\mR^n,\mR^k]_{\mathrm{prop}}$ is in bijection with $\pi_{n-1}(S^{k-1})$. The set $ [\mR^n,\mR]_{\mathrm{prop}}$ has two elements if $n>1$ and four elements if $n=1$. 
\end{corollary}

A proper map between non-compact and locally compact Hausdorff spaces extends to a continuous map between the one point compactifications by sending infinity to infinity. Similarly a homotopy of proper maps induces a homotopy in the one-point compactification. 

The one point compactification of a real vector bundle $E\rightarrow M$ over a compact Hausdorff space $M$ equals the Thom space $\Th(E)$ of the vector bundle and the one point compactification of $\mR^k$ is homeomorphic to $S^k$ by stereographic projection. Thus we obtain a map $Q:[E,\mR^k]_{\mathrm{prop}}\rightarrow \langle\Th(E),S^k\rangle$. In Section~\ref{sec:suspensions} we show that the map $Q$ is bijective in a range. If $E=\mR^n$ the map $Q$ is nothing but the suspension $\pi_{n-1}(S^{k-1})\rightarrow \pi_n(S^k)$ under the identification of $[\mR^n,\mR^k]_{\mathrm{prop}}$ and $\pi_{n-1}(S^{k-1})$ of Corollary~\ref{cor:sphere}.

For $l$ sufficiently the sets $[E\oplus \mR^l,\mR^{k+l}]_{\mathrm{prop}}$ and $[E\oplus \mR^{l+1},\mR^{k+l+1}]$ are in bijection. Thus it makes sense to define the stable proper homotopy classes as $[E,\mR^k]^S_{\mathrm{prop}}=\lim_{l\rightarrow \infty}[E\oplus \mR^l,\mR^{k+l}]$, which are in bijection with the stable cohomotopy groups $\pi^k_S(\Th(E))$, cf.~Corollary~\ref{cor:stable}. Using Atiyah duality we obtain the following result. 

\begin{theorem}
\label{thm:atiyah}
Let $M$ be an $m$-dimensional compact manifold with boundary $\partial M$. Let $E$ be the normal bundle of $M$ of some embedding of $M$ into $\mR^{m+n}$. Then there is a bijection of $[E,\mR^k]^S_{\mathrm{prop}}$ with $\pi_{n+m-k}^S(M/\partial M)$.
\end{theorem}

In Section~\ref{sec:speculation} we speculate on the classification problem in the case that $E$ is an arbitrary open manifold. 

\subsection{Acknowledgements} I would like to thank Alberto Abbondandolo, Hansj\"org Geiges, Gijs Heuts, and Federica Pasquotto for discussions on the content of this paper. This research was supported by NWA startimpuls - 400.17.608.
\section{The proof of Theorem~\ref{maintheorem}}

For the remainder of the paper $E\rightarrow M$ denotes a normed real vector bundle over a compact Hausdorff space $M$. The associated sphere and disc bundles of radius $R>0$ are
$$S_R(E)=\{v\in E\,|\,\norm{v}=R\}\, \quad\text{and}\quad B_R(E)=\{v\in E\,|\,\norm{v}<R\}.$$ 
We write $S(E)$ for $S_1(E)$ and $B(E)$ for $B_1(E)$. 

Given a homotopy $F:[0,1]\times S(E)\rightarrow S^{k-1}$ we define the homotopy $PF:[0,1]\times E\rightarrow \mR^k$ of proper maps by
$$
PF(t,v)=\begin{cases}
\norm{v}F\left(t,\frac{v}{\norm{v}}\right)\qquad &v\not=0\\
0\qquad &v=0. 
\end{cases}
$$
Compact subsets of $E$ are characterized as follows: A subset $K\subseteq E$ is compact if and only if it is closed and bounded. Here bounded means that $K\subseteq B_R(E)$ for some $R>0$. Let us prove that $PF$ is proper. Let $C\subseteq \mR^k$ be compact. Then $C$ and hence $PF^{-1}(C)$ are closed as $PF$ is continuous. Compact subsets of $\mR^k$ are the closed and bounded subsets, hence the set $C$ is contained in $B_r(\mR^k)$ for some $r>0$. Thus $PF^{-1}(C)$ is a closed subset contained the bounded set $[0,1]\times B_r(E)$ hence is compact. We conclude that that $PF$ is proper. 

The same construction assigns to a map $f:S(E)\rightarrow S^{k-1}$ a proper map $Pf:E\rightarrow \mR^k$ and  it therefore induces a map $P:[S(E),S^{k-1}]\rightarrow [E,\mR^k]_{\mathrm{prop}}$. 

We will show that $P$ is bijective. Let us start with the injectivity. We need to show that if $g_0=Pf_0$ and $g_1=Pf_1$ are homotopic as proper maps, that this implies that $f_0$ and $f_1$ are homotopic. Let $G:[0,1]\times E\rightarrow \mR^k$ be a homotopy  of proper maps between $g_0$ and $g_1$. Then for any $r>0$ there exists an $R>0$ such that $G^{-1}(B_r(\mR^k))\subseteq [0,1]\times B_R(E)$. It follows that for any $v\in E$ with $\norm{v}=R$ that $G(t,v)\not=0$. The map $F:[0,1]\times S(E)\rightarrow S^{k-1}$ given by
$$
F(t,x)=\frac{G(t,Rx)}{\norm{G(t,Rx)}},
$$
is a homotopy between $f_0$ and $f_1$ hence $P$ is injective.

To show that $P$ is surjective, we need to show that, given a proper map $g:E\rightarrow \mR^k$, there exists a homotopy of proper maps from $g$ to $Pf$, where $f$ is some map $f:S(E)\rightarrow S^{k-1}$. As $g$ is proper, there exists an $R>0$ such that $g^{-1}(B(\mR^k))\subseteq B_R(E)$. The sphere bundle $S_R(E)$ is compact, hence there exists an $r\geq 1$ such that
$$
1\leq\norm{g(v)}\leq r,\quad\text{for all}\quad v\in S_R(E). 
$$
Consider the map $h:\mR^k\rightarrow \mR^k$ defined by
$$
h(x)=\begin{cases}x\qquad &\norm{x}\leq 1\\
\frac{x}{\norm{x}}\qquad &1\leq \norm{x}\leq r\\
\frac{x}{r}\qquad &r\leq \norm{x},
\end{cases}
$$
and define $g_1:E\rightarrow \mR^k$ by $g_1(v)=h\circ g(Rv)$. As $h$ is homotopic as proper map to the identity via $(t,x)\mapsto (1-t) h(x)+x$ it follows that he map $g_1$ is proper homotopic to $g$. Note that $g_1(B(E))\subseteq B(\mR^k)$ and $g_1(S(E))\subseteq S^{k-1}$. Define $G_1:[0,1]\times E\rightarrow E$ by
$$
G_1(t,v)=\begin{cases}\frac{\norm{g_1(v)}}{\norm{g_1\left((1-t) v+\frac{t}{\norm{v}}v\right)}}g_1\left((1-t)v+\frac{t}{\norm{v}}v\right) \qquad &\norm v\geq 1\\
g_1(v)\qquad &\norm v\leq 1.
\end{cases}
$$
The equation $\norm{G_1(t,v)}=\norm{g_1(v)}$ implies that $G_1^{-1}(B_s(\mR^k))\subseteq [0,1]\times g_1^{-1}(B_s(\mR^k))$ for all $s$ and hence that $G_1$ is proper if $g_1$ is. Thus $g_2:E\rightarrow \mR^k$ given by $g_2(x)=G_1(1,x)$ is a proper map that is proper homotopic to $g$. Let $f:S(E)\rightarrow S^{k-1}$ be the map obtained by restriction of $g_2$. Consider 
$$
G_2(t,v)=(1-t)g_2(v)+t Pf(v). 
$$
We want to prove that $G_2$ is proper. Note that $G_2$ sends $\overline{ [0,1]\times B(E)}$ to $\overline{B(\mR^k)}$ and $[0,1]\times (E\setminus B(E))$ to $\mR^k\setminus B(\mR^k)$. The map $G_2\bigr|_{\overline{[0,1]\times B(E)}}: \overline{[0,1]\times B(E)}\rightarrow \overline{B(\mR^k)}$ is proper as the domain is compact. We conclude that $G_2$ is proper if and only if $$G_2\bigr|_{[0,1]\times (E\setminus B(E))}: [0,1]\times (E\setminus B(E))\rightarrow \mR^k\setminus B(\mR^k)$$ is proper.  

For this it is sufficient to show that for all $s>1$ there exists an $S>1$ such that $G_2^{-1}(B_s(\mR^k)\setminus B(\mR^k))\subseteq [0,1]\times (B_S(E)\setminus B(E))$. Note that for all $(t,v)\in [0,1]\times (E\setminus B(E))$ we have that
$$
G_2(t,v)=\left((1-t)\norm{g_1(v)}+t \norm v\right) f\left(\frac{v}{\norm v}\right).
$$
Consider all $(t,v)\in [0,1]\times (E\setminus B(E))$ such that $\norm{G_2(t,v)}\leq s.$ As $\norm {f(\frac{v}{\norm v})}=1$ this amounts to 
$$
(1-t)\norm {g_1(v)}+t\norm{v}\leq s.
$$
Suppose the set of solutions of this equation is not contained in $[0,1]\times (B_S(E)\setminus B(E))$ for any $S$. Then we have a sequence $(t_n,v_n)$ of solutions such that $\norm{v_n}\geq n$. Without loss of generality we take a subsequence such that $t_n$ converges to $t$ by the compactness of $[0,1]$. This subsequence will also satisfy $\norm{v_n}\geq n$. If $t>0$, then there exists an $N$ such that for all $n\geq N$ we have $t_n>\frac{t}{2}$ and 
$$
\norm{v_n}\leq \frac{1}{t_n}\left((1-t_n)\norm {g_1(v_n)}+t_n\norm{v_n}\right)\leq \frac{2s}{t},
$$
which contradicts the unboundedness of $v_n$. If $t=0$ then there exists an $N$ such that for all $n\geq N $ the sequence satisfies $t_n<\frac{1}{2}$ and
$$
\norm {g_1(v_n)}\leq\frac{1}{(1-t_n)}\left((1-t_n)\norm {g_1(v_n)}+t_n\norm{v_n}\right)\leq 2s.
$$
The sequence $g_1(v_n)$ is therefore bounded and as the map $g_1$ is proper it follows that the sequence $v_n$ is also bounded. This contradicts the assumption that $v_n$ is unbounded. This means that $G_2\bigr|_{[0,1]\times (E\setminus B(E))}$ is proper. Thus $Pf$ is proper homotopic to $g$ and $P:[S(E),S^{k-1}]\rightarrow [E,\mR^k]_{\mathrm{prop}}$ is surjective. We have already shown that $P$ is injective and Theorem~\ref{maintheorem} follows. 

\section{The one point compactification and stable (co)homotopy}
\label{sec:suspensions}

Recall that the one point compactification of a non-compact, locally compact Hausdorff space $X$ is the space $X^*=X\cup\{\infty\}$ equipped with the following topology. All open sets $U$ of $X$ are declared open in $X^*$ along with all sets of the form $(X\setminus C)\cup \{\infty\}$ for all compact sets $C$ in $X$. Proper maps between non-compact, locally compact Hausdorff spaces induce continuous maps between the one point compactifications by imposing that $\infty$ is mapped to $\infty$. A homotopy of proper maps $F:[0,1]\times X\rightarrow Y$ induces a continuous map $F^*:([0,1]\times X)^*\rightarrow Y^*$. But $([0,1]\times X)^*\cong [0,1]\times X^*/[0,1]\times \{\infty\}$. By the universal property of the quotient topology we therefore also obtain a continuous map $F^*:[0,1]\times X^*\rightarrow Y^*$ which sends every $(t,\infty)$ to $\infty$. A homotopy of proper maps between unbased spaces is mapped to a based homotopy between the based spaces. The one point compactification of a vector bundle $E\rightarrow M$ over a compact Hausdorff space $M$ is the Thom space of the bundle and we will write $\Th(E):=E^*$. Stereographic projection shows that $(\mR^k)^*\cong S^k$ and in more generality it holds that $\Th(E)\cong \overline{B(E)}/S(E)$. As was mentioned in the introduction, the map that forgets the basepoint induces a bijection between $\langle \Th(E),S^k\rangle$ and $[\Th(E),S^k]$ if $k\geq 1$. Thus from the one point compactification we obtain a map $Q:[E,\mR^k]_{\mathrm {prop}}\rightarrow [\Th(E),S^k]$. 

In our setting there are three suspension maps, which we all denote by $S$. To a map $g:S(E)\rightarrow S^{k-1}$ we associate the map $S f:S(E\oplus \mR)\rightarrow S^k$ by
$$
Sg(v,s)=\begin{cases}
(\norm v g\left(\frac{v}{\norm v}\right),s)\qquad &v\not=0\\
(0,s)&v=0.
\end{cases}
$$ 
To a map $f:\Th(E)\rightarrow S^k$ we associate the map $Sf:\Th(E\oplus \mR)\rightarrow S^{k+1}$  via the same formula and to a proper map $f:E\rightarrow \mR^k$ we associate the proper map $Sf:E\oplus \mR\rightarrow \mR^{k+1}$, by $Sf(x,s)=(f(x),s)$. The following diagram is commutative
\begin{equation}
\label{diag}
\begin{gathered}
\xymatrix{
[S(E),S^{k-1}]\ar[d]_S\ar[r]^P&[E,\mR^k]_{\mathrm{prop}}\ar[d]_S\ar[r]^Q&[\Th(E),S^k]\ar[d]_S\\
[S(E\oplus\mR),S^{k}]\ar[r]^P&[E\oplus\mR,\mR^{k+1}]_{\mathrm{prop}}\ar[r]^Q&[\Th(E\oplus\mR),S^{k+1}].
}
\end{gathered}
\end{equation}
In the proof of Theorem~\ref{maintheorem} we saw that the maps $P$ are bijections. We wonder when the other maps in the diagram are bijective. 

Let us now assume that $k\geq 2$, that $M$ is a finite connected CW-complex of dimension $m$, and that $E$ is a vector bundle of rank $n$. Since $k\geq 2$ based and unbased homotopy classes into $S^{k-1}$ and $S^k$ coincide, as well as based or unbased proper homotopy classes into $\mR^k$ and $\mR^{k+1}$. We denote by $$\pi^{k-1}(S(E)):=\langle S(E),S^{k-1}\rangle$$ the $(k-1)$-th cohomotopy set of $S(E)$. We refer to \cite[Chapter VII]{Hu:1959ve} for information on the cohomotopy sets we use below. The cohomotopy set $\pi^{k-1}(S(E))$ is not always a group, but only if $m+n\leq 2k-3$. We investigate the long exact sequence of the pair $(\overline{B(E)},S(E))$ if $m+n\leq 2k-3$.
\begin{equation}
\label{eq:longexact}
\pi^{k-1}(\overline{B(E)})\rightarrow \pi^{k-1}(S(E))\xrightarrow{\delta} \pi^k(\overline{B(E)},S(E))\rightarrow \pi^k(\overline{B(E)})
\end{equation}
Since $\overline{B(E)}$ deformation retracts to $M$ and $S^{k-1}$ is $(k-2)$-connected, we see that if $m\leq k-2$ there are  isomorphisms 
$$
\pi^{k-1}(\overline{B(E)})\cong\pi^{k-1}(M)\cong 0\quad\text{ and}\quad \pi^{k}(\overline{B(E)})\cong\pi^{k}(M)\cong 0.
$$ 
Thus we conclude that for $2k\geq m+3+\max(n,m+1)$ there is an isomorphism $\pi^{k-1}(S(E))\cong\pi^k(\overline{B(E)},S(E))$. The relative cohomotopy set is the cohomotopy set of the quotient for nice spaces, thus $\pi^{k}(\overline{B(E)},S(E))=\pi^k(\overline{B(E)}/S(E))=\pi^k(\Th(E))$. The coboundary map is an isomorphism $\pi^{k-1}(S(E))\cong\pi^k(\Th(E))$ in the dimension range. Let us consider the based version of Diagram~\eqref{diag}
\begin{equation}
\begin{gathered}
\xymatrix{
\pi^{k-1}(S(E)) \ar[d]_S\ar[r]^{QP}&\pi^{k}(\Th(E))\ar[d]_S\\
\pi^{k}(S(E\oplus \mR))\ar[r]^{QP}&\pi^{k+1}(\Th(E\oplus R)).
}
\end{gathered}
\end{equation}
The horizontal maps can be identified with the coboundary map $\delta$ in~\eqref{eq:longexact} and therefore the horizontal maps are isomorphisms in the right dimension range. Freudenthal's suspension Theorem, cf.~\cite{Kochman:1996ir}, states that if $m+n\leq 2k-2$, that the suspension map $\pi^k(\Th(E))\rightarrow \pi^{k+1}(S\Th(E))\cong\pi^{k+1}(\Th(E\oplus\mR))$ is an isomorphism. Combining all information we have now gives us
\begin{theorem}
\label{thm:iso}
Let $M$ be a finite $CW$ complex of dimension $m$ and $E$ a vector bundle over $M$ of rank $n$. Let $k\geq 2$ and suppose that $2k\geq m+3+\max(n,m+1)$. Then all maps in Diagram~\eqref{diag} are bijections. 
\end{theorem}

This theorem expresses a stability phenomenon: For all $l$ sufficiently large the map $Q$ induces bijections $[E\oplus \mR^l,\mR^{k+l}]_{\mathrm{prop}}\rightarrow [E\oplus \mR^{l+1},\mR^{k+l+1}]_{\mathrm{prop}}$. We define the stable homotopy classes of proper maps as
$$
[E,\mathrm{\mR^k}]^S_{\mathrm{prop}}=\lim_{l\rightarrow \infty}[E\oplus \mR^l,\mR^{k+l}]_{\mathrm{prop}}.
$$
Recall that the stable homotopy and cohomotopy groups of a space $X$ are similarly defined
$$
\pi_k^S(X)=\lim_{l\rightarrow \infty}\langle S^l S^k,S^l X\rangle\quad \text{and}\quad \pi_S^k(X)=\lim_{l\rightarrow \infty}\langle S^l X,S^l S^k\rangle.
$$
A direct corollary of Theorem~\ref{thm:iso} is then
 \begin{corollary}
 \label{cor:stable}
 Let $M$ be a finite $CW$ complex of dimension $m$ and $E$ a vector bundle of rank $n$. Then $Q$ induces a bijection $[E,\mR^k]^S_{\mathrm{prop}}$ with $\pi^k_S(\Th(E))$
 \end{corollary}

Stable homotopy and cohomotopy groups are related via Spanier-Whitehead duality, which we recall now. We refer to the original references~\cite{Spanier:1956hv,Spanier:1955hg} for these statements. Let $i:X\rightarrow S^N$ be a sufficiently nice embedding of a sufficiently nice space $X$ into a sphere (e.g. a smooth embedding of a compact manifold, or the inclusion of a CW complex as a subcomplex). Then the space $D_NX=S^N\setminus i(X)$ is a Spanier-Whitehead dual of $X$. The stable homotopy type of $D_NX$ is well defined: It is independent of the dimension $N$ and the choice of embedding. The fundamental result is that $\lim_{l\rightarrow \infty}[S^lX,S^lY]$ is in bijection with $\lim_{l\rightarrow \infty}[S^lD_NY,S^lD_NX]$. In particular, the stable cohomotopy groups of $X$ are the stable homotopy groups of $D_NX$ with a dimension shift. Now let us assume that $M$ is a compact manifold with boundary $\partial M$. There is a unique (up to isotopy) embedding of $M$ into $\mR^{m+n}$ for $n$ sufficiently large. Let $E$ be the normal bundle of such an embedding, i.e. let $E$ be the stable normal bundle of $M$. Atiyah~\cite[Proposition 3.2]{Atiyah:1961dp} showed that $SD_{m+n}(M/\partial M)\simeq\Th(E)$. If the boundary $\partial M$ is empty we should interpret $M/\partial M$ as $M$ with a disjoint basepoint added. The Spanier-Whitehead dual of a sphere is $D_{n+m} S^{n+m-k}=S^{k-1}$. We have
\begin{align*}
\pi_{n+m-k}^S(M/\partial M):&=\lim_{l\rightarrow\infty}\langle S^lS^{n+m-k},S^l M/\partial M\rangle\\
&=\lim_{l\rightarrow \infty}\langle S^l D_{n+m}(M/\partial M),S^l D_{n+m} S^{n+m-k}\rangle\\
&=\lim_{l\rightarrow\infty}\langle S^{l-1}\Th(E),S^l S^{k-1}\rangle\\
&=\lim_{l\rightarrow \infty}\langle S^l\Th(E),S^l S^k\rangle\\
&=\pi^k_S(\Th(E)).
\end{align*}
Theorem~\ref{thm:iso} states that $\pi^k_S(\Th(E))$ is in bijection with $[E,\mR^k]^S_{\mathrm{prop}}$. We have proven Theorem~\ref{thm:atiyah}.

\section{Framed submanifolds and cobordisms}
\label{framed}
Pontryagin~\cite{Pontryagin:1959vj} showed that homotopy classes of maps $M\rightarrow S^k$, where $M$ is a closed manifold, are in one to one correspondence with framed cobordism classes of $(n-k)$-dimensional manifolds in $M$. Framed cobordism classes are also invariants of homotopy classes of proper maps $E\rightarrow \mR^k$ but they are not complete, nor is every cobordism classed realized by some proper map. In this section we discuss this. 

\subsection{Invariants of proper maps: framed submanifolds and cobordisms}
Let $M$ be a smooth $m$-dimensional manifold and $N$ be a connected oriented smooth $k$-dimensional manifold. Every continuous proper map is homotopic as a proper map to a smooth proper map, hence we can consider only proper smooth maps in the proper homotopy classification of proper maps. Suppose that $f:M\rightarrow N$ is a smooth map that is proper. Proper maps between manifolds are closed maps. The set of critical points is closed, hence the set of regular values of a proper map is open. By Sard/Brown's theorem regular values of $f$ are residual, and by Baire's category theorem it follows that the regular values are dense. Thus the set of regular values of a proper map is open and dense. An example of smooth function $\mR_{>0}\rightarrow \mR$ whose set of regular values is not open is $x\mapsto \frac{1}{x(2+\sin(x))}$, but of course this map fails to be proper and closed. The preimage of a regular value $y$ is a closed submanifold $X=f^{-1}(y)$ of dimension $m-k$. Such a manifold can be \emph{framed}: Let $e_1,\ldots,e_k$ be a basis of $T_yN$ that is compatible with the orientation of $N$. Then for every $x\in X$ the differential of $f$ induces an isomorphism $df_x:N_xX\rightarrow T_yN$ of the normal space $N_xX$ to $X$ at $x$ with $T_yN$. Then $(\nu_f)_x=((df_x)^{-1}(e_1),\ldots, (df_x)^{-1}(e_k))$ is an ordered basis of the normal space $N_xX$ at $x$. Letting $x$ vary, this patches together to a map $\nu_f$ that trivializes the normal bundle of $X$. The map $\nu_f$ is called the framing of $X$. We call $(X,\nu_f)$ a Pontryagin manifold of $f$ and it depends on the choices we made.  

Let $F:[0,1]\times M\rightarrow N$ be a homotopy of proper maps between $f_0=F(0,\cdot)$ and $f_1=F_1(1,\cdot)$. By a reparametrization of the homotopy variable, we may assume that $F(t,x)=f_0(x)$ and $F(1-t,x)=f_1(x)$ for $t$ small. If $y$ is a regular value of the maps $F,f_0$ and $f_1$ simultaneously then $(W=F^{-1}(y),\nu_F)$ is a framed compact submanifold with framed boundary $(X_0=f_0^{-1}(y),\nu_{f_0})$ and $(X_1=f_0^{-1}(y),\nu_{f_1})$. The framed manifold $(W,\nu_F)$ is a framed cobordism between the framed manifolds $(X_0,\nu_{f_0})$ and $(X_1,\nu_{f_1})$. 
Being framed cobordant defines an equivalence relation on the set of framed submanifolds and the framed cobordism class of a Pontryagin manifold of a proper map $f:M\rightarrow N$ does not depend on the choice of regular value $y$ and the choice of oriented basis of $T_yN$ and is an invariant of the proper homotopy class of $f$. We denote the set of framed $(m-k)$-dimensional closed submanifolds of $M$ up to framed cobordism by $\frc_{m-k}(M)$. 

\subsection{The Pontryagin-Thom construction}

The framed cobordism class of the preimage of a regular value is in some cases enough to recover the homotopy class of the map: Suppose $M$ is closed and $(X,\nu)$ is a $(m-k)$-dimensional framed submanifold. Out of this data we can construct a (proper) map $f:M\rightarrow S^k$, for which $(X,\nu)$ is a Pontryagin manifold: We define $f$ to map $X$ to the northpole $y$ of $S^k$ and describe what happens in a tubular neighborhood of $X$. The framing $\nu$ defines, for each point $x\in X$ a diffeomorphism of the normal space around $x$ to a neighborhood of $y$. We use this to extend the map to the tubular neighborhood of $X$ in $M$. One can arrange this in such a way that if one approaches the boundary of the tubular neighborhood, the image under $f$ converges to the south pole. The map $f$ can now be extended to the whole of $M$ by mapping everything outside the tubular neighborhood to the south pole. The northpole is a regular value for $f$ and the Pontryagin manifold at the north pole is exactly the framed manifold $(X,\nu)$. This construction also works for framed cobordisms. This proves the following theorem.

\begin{theorem}
The Pontryagin-Thom construction gives a one to one correspondence between the set $\frc_{m-k}(M)$ of framed cobordisms in $M$ and the set $ [M,S^k]$ of homotopy classes of maps from $M$ to $S^k$. 
\end{theorem}

For more details of the Pontryagin-Thom construction in this classical setting we refer to Milnor~\cite{Milnor:_FIovTN9} and Pontryagin~\cite{Pontryagin:1959vj}.

\subsection{How good of an invariant is the framed cobordism class of the Pontryagin manifold of a proper map?}
\label{sec:counterexample}
A proper map $f:\mR^n\rightarrow \mR^k$ is proper homotopic to a map $Pg$, where $g:S^{n-1}\rightarrow S^{k-1}$. Let $y\in S^{n-1}$ be a regular value of $g$, and $(X,\nu_g)$ be the Pontryagin manifold of $g$ at $y$. The value $y\in S^{k-1}\subseteq \mR^k$ is also a regular for the map $Pg$. The Pontryagin manifold of $Pg$ at $x$ is $(X,(\nu_{g},\mu))$, where $\mu$ is the last component of the framing which points radially outward from the sphere. So a framed submanifold cannot occur as a Pontryagin manifold if it is not framed cobordant to a framed submanifold that lies on a sphere where the last component of the framing is radially pointing outward.

Let us discuss an explicit example of a framed manifold that does not occur as the Pontryagin manifold of a map. Consider the submanifold $X=\{-1,1\}\subseteq \mR$ with framing $\nu_{-1}=\nu_1=\frac{\partial}{\partial t}$. Then $(X,\nu)$ cannot occur as the preimage of a regular value of a proper map $f:\mR\rightarrow \mR$. Suppose on the contrary that such a map exists with $f(-1)=f(1)$. From the framing and Taylor's theorem we see that there exists an $\epsilon>0$ such that $f(-1+\epsilon)>f(-1)$ and $f(1-\epsilon)<f(1)$. The intermediate value theorem then gives the existence of another point $-1+\epsilon<p<1-\epsilon$ such that $f(p)=f(-1)=f(1)$. Hence $f^{-1}(f(1))\not=X$ and we conclude that there does not exist a proper $f:\mR\rightarrow \mR$ with $(X,\nu)$ as a Pontryagin manifold. 

But there are also framed submanifolds which are framed cobordant to a framed submanifold which is contained in the unit sphere and has a framing with last component pointing radially outward which do not arise from proper maps. To see this, consider the manifold $X=\{-2,-1,1,2\}$ with framing $\nu(-2)=\nu(2)=\frac{\partial}{\partial t}$ and $\nu(-1)=\nu(1)=-\frac{\partial}{\partial t}$. Then $(X,\nu)$ is framed cobordant to the empty set, however it cannot occur as the Pontryagin manifold of a proper map: If $y$ is the regular value for which $(X,\nu)$ is hypothetically the Pontryagin manifold at $y$, there must be a point $x \in (-1,1)$ such that $f(x)=y$, by the same reasoning as above. 

Finally we discuss the fact that the invariant is not complete. The maps $f,g:\mR\rightarrow \mR$ given by $f(x)=x^2$ and $g(x)=-x^2$ are not proper homotopic. However as the maps are not surjective the framed cobordism class of both maps is the empty manifold. Hence the framed cobordism class of a regular value cannot distinguish these maps. Here is a more complicated example: Let $f,g:S^3\rightarrow S^2$ be the Hopf map and the Hopf map precomposed with a degree $-1$ map of $S^3$ respectively. These maps represent $+1$ and $-1$ in $\pi_3(S^2)\cong \mZ$ and are not homotopic. By Theorem~\ref{maintheorem} $Pf$ and $Pg$ are therefore not proper homotopic, however their Pontryagin manifolds are framed cobordant. To see this note that $QPf$ and $QPg$ are the suspensions of $f$ and $g$ and the Pontryagin manifolds of $QPf$ and $QPg$ can be identified with those of $Pf$ and $Pg$. The suspension map $S:\pi_3(S^2)\cong \mZ\rightarrow \pi_4(S^3)\cong \mZ/2\mZ$ is the  reduction modulo $2$. The maps $QPf$ and $QPg$ are homotopic, so their Pontryagin manifolds must be framed cobordant. But this implies that the Pontryagin manifolds of $PF$ and $PG$ are framed cobordant by general position. 

\subsection{Open manifolds and the Pontryagin-Thom construction}

\label{sec:speculation}
In Corollary~\ref{cor:stable} we have seen that the homotopy classes of proper maps out of vector bundles stabilizes. We expect that if $M$ is an arbitrary open manifold the homotopy classes of proper maps $[M\times \mR^l,\mR^{k+l}]_{\mathrm{prop}}$ stabilizes when $l\rightarrow \infty$. This suggests that there is a stable Pontryagin-Thom construction for proper maps. A framed submanifold $(X,\nu)\in \frc_{m-k}(M)$ gives rise to a framed submanifold $((X,0),\nu\oplus \mu)\in \frc_{m-k}(M\times \mR^l)$ via stabilization. The framing $\nu\oplus \mu $ extends the framing $\nu$ with a fixed basis $\mu$ of $\mR^l$. Since we expect that the homotopy classes of proper maps stabilizes we also expect that there is a well defined stable bijective Pontryagin-Thom construction $\frc_{m-k}(M\times \mR^l)\rightarrow [M\times \mR^l,\mR^{l+k}]_{\mathrm{prop}}$ for $l$ sufficiently large.

\bibliographystyle{abbrv} \bibliography{biblio}

\begin{thebibliography}{10}

\bibitem{Abbondandolo:2018fb}
A.~Abbondandolo and T.~O. Rot.
\newblock {On the homotopy classification of proper Fredholm maps into a
  Hilbert space}.
\newblock {\em Journal f{\"u}r die Reine und Angewandte Mathematik}, (Ahead of
  Print) 2018.

\bibitem{Atiyah:1961dp}
M.~F. Atiyah.
\newblock {Thom complexes}.
\newblock {\em Proceedings of the London Mathematical Society. Third Series},
  11(1):291--310, 1961.

\bibitem{Elworthy:1970vp}
K.~D. Elworthy and A.~J. Tromba.
\newblock {Differential structures and Fredholm maps on Banach manifolds}.
\newblock In {\em Global Analysis (Proc. Sympos. Pure Math., Vol. XV, Berkeley,
  Calif., 1968)}, pages 45--94. Amer. Math. Soc., Providence, R.I., 1970.

\bibitem{Geba:1969df}
K.~G{\k{e}}ba.
\newblock {Fredholm $\sigma$ -proper maps of Banach spaces}.
\newblock {\em Polska Akademia Nauk. Fundamenta Mathematicae}, 64(3):341--373,
  1969.

\bibitem{Geiges:2018fo}
H.~Geiges.
\newblock {Isotopies vis-\`a-vis level-preserving embeddings}.
\newblock {\em Archiv der Mathematik}, 110(2):197--200, 2018.

\bibitem{Hatcher_Topology}
A.~Hatcher.
\newblock {\em {Algebraic topology}}.
\newblock Cambridge University Press, Cambridge, 2002.

\bibitem{Hu:1959ve}
S.~T. Hu.
\newblock {\em {Homotopy theory}}.
\newblock Pure and Applied Mathematics, Vol. VIII. Academic Press, New
  York-London, 1959.

\bibitem{Kochman:1996ir}
S.~O. Kochman.
\newblock {\em {Bordism, stable homotopy and Adams spectral sequences}},
  volume~7 of {\em Fields Institute Monographs}.
\newblock American Mathematical Society, Providence, RI, Providence, Rhode
  Island, 1996.

\bibitem{Milnor:_FIovTN9}
J.~W. Milnor.
\newblock {Topology from the differentiable viewpoint}.
\newblock pages ix+65, 1965.

\bibitem{Pontryagin:1959vj}
L.~S. Pontryagin.
\newblock {Smooth manifolds and their applications in homotopy theory}.
\newblock In {\em American Mathematical Society Translations, Ser. 2, Vol. 11},
  pages 1--114. American Mathematical Society, Providence, R.I., 1959.

\bibitem{Spanier:1956hv}
E.~H. Spanier.
\newblock {Duality and S-theory}.
\newblock {\em Bulletin of the American Mathematical Society}, 62(3):194--203,
  1956.

\bibitem{Spanier:1955hg}
E.~H. Spanier and J.~H.~C. Whitehead.
\newblock {Duality in homotopy theory}.
\newblock {\em Mathematika. A Journal of Pure and Applied Mathematics},
  2:56--80, 1955.

\end{thebibliography}
\end{document}